\documentclass[journal]{IEEEtran}

\IEEEoverridecommandlockouts                              % This command is only needed if 
\usepackage{amsmath}
\usepackage{bbm}

\usepackage{nomencl}

\usepackage[fleqn]{nccmath}
\usepackage{graphicx}
\usepackage{amsfonts}

\usepackage{array}
\usepackage{booktabs}
\usepackage{makecell}
\usepackage[utf8]{inputenc}
\usepackage{blkarray}

\usepackage{caption}
\usepackage{subcaption}
\usepackage[hyphens,spaces]{url}
\usepackage{tikz}
\usetikzlibrary{shapes, arrows, chains}
\usepackage[flushleft]{threeparttable}

\usepackage[ruled,vlined]{algorithm2e}
\usepackage[flushleft]{threeparttable}
\usepackage{algpseudocode}

\title{
Joint Chance-constrained Game for Coordinating Renewable Microgrids with Service Delivery Risk: A Bayesian Optimization Approach}

\author{Yifu Ding, Benjamin Hobbs% <-this % stops a space
\thanks{$^{1}$Yifu Ding is with the Massachusetts Institute of Technology Energy Initiative, USA.
        {\tt\small yifuding@mit.edu}}%
\thanks{$^{2}$Benjamin Hobbs is with Whiting School of Engineering, Johns Hopkins University, USA
        {\tt\small bhobbs@jhu.edu}}%
}

\begin{document}

\maketitle
\thispagestyle{empty}
\pagestyle{empty}

%%%%%%%%%%%%%%%%%%%%%%%%%%%%%%%%%%%%%%%%%%%%%%%%%%%%%%%%%%%%%%%%%%%%%%%%%%%%%%%%
\begin{abstract}

Microgrids incorporate distributed energy resources (DERs) and flexible loads, which can bid for energy and reserve service provision in competitive auctions. However, due to uncertain renewable generations such as solar power, microgrids might fail to deliver reserve services and breach day-ahead contracts in real time. This undermines system reliability and leads to grid contingencies. This paper presents a distributionally robust joint chance-constrained (DRJCC) game framework that facilitates the coordination of renewable microgrids in the bilateral market, ensuring reliable energy and reserve service provisions. Leveraging historical forecast error samples, the reserve bidding strategy of each microgrid is formulated into a two-stage Wasserstein-metrics distribution robust optimization (DRO) model. A JCC regulates the under-delivered reserve capacity of all microgrids in a non-cooperative game. Considering the unspecific correlation among players, an end-to-end Bayesian optimization method swiftly approximates the individual contract violation rates of microgrids, given the unique game equilibria. The proposed Bayesian optimization game framework is simulated with up to 14 players using the California power market data. Results show the proposed method can effectively regulate the under-delivered reserve of all microgrids jointly and secure the profit of each microgrid in the reserve market.


\end{abstract}

\IEEEpeerreviewmaketitle

\section{Introduction}

The transition to the net-zero energy system simulates the adoption of behind-the-meter generations such as solar PV panels \cite{department_for_business_energy__industrial_strategy_powering_nodate}. These assets in the distribution network can provide flexibility but cannot directly receive price signals from the wholesale power or ancillary markets. An organized ancillary services market is usually subject to regulation or bilateral contracts with the system operator, sometimes involving an auction procedure \cite{keay_electricity_2013}. These distributed energy resources thus could be aggregated as a microgrid to reach the minimum capacity for market participation. Considering they have no explicit energy and power capacity, the performance assessment must be conducted to meet the technical requirements for service provisions and reach a binding contract. Nevertheless, weather-dependent renewable generations and stochastic demands of microgrids might lead to service under-delivery and breach contracts, undermining system reliability and leading to grid contingencies. 

Game-theoretical frameworks usually simulate a competitive procedure such as auctions, where market participators have their objectives. Their biddings depend on risk appetites and the market environment until market equilibrium is reached. This market dynamics can be modeled mainly by either of two methods. They are the iterative or heuristic algorithm to search for the equilibrium, such as the best response algorithm, and mathematical programming with equilibrium constraints (MPEC) \cite{gabriel_equilibria_2013}. 

The former method simulates each player's strategy iteratively and searches for a market equilibrium until all players' solutions converge to a near-optimal point \cite{delatorre_finding_2004, li_data-driven_2022}. However, this best-response method is a heuristic method and only guarantees convergence to an efficient equilibrium for certain utility functions. The game result and computation time are also sensitive to the starting point \cite{saad_game-theoretic_2012}. The MPEC model based on the Karush-Kuhn-Tucker (KKT) rule secures a fixed game equilibrium when at least a feasible solution exists. Still, it results in the non-convex optimization problem requiring additional mathematical reformulations \cite{saguan_market_2006}.

The power markets in the U.S. and Europe are usually two-settlement or multi-settlement. The two-settlement market structures consist of the day-ahead forward and real-time markets. The intra-day stage will also be included if it is a multi-settlement market. In the literature, these markets are simulated in bilevel or trilevel models involving independent system operators (ISOs), distribution system operators (DSOs), and service providers \cite{jiang_flexibility_2022}, \cite{hobbs_strategic_2000-1} - \cite{li_distributed_2023}. In competitive market models, the interactions of the ISO and DSO (or microgrids in this paper) are often modeled by the Stackelberg and Nash bargaining games \cite{jiang_flexibility_2022, hobbs_strategic_2000-1, bruninx_interaction_2020}, while the market competition of DSOs is modeled by the Nash-Cournot games \cite{li_risk-averse_2017, li_distributed_2023}. 

%\textcolor{black}{The Stackelberg game models DSOs as relatively small consumers competing in the market and passively responding to ISO’s action \cite{jiang_flexibility_2022}. In contrast, the Nash bargaining game could model the large DSOs who can negotiate prices and cooperate with ISO—for example, ref. \cite{bruninx_interaction_2020} models the interactions of aggregators and demand response providers. This work shows that the existing outcome of a Stackelberg Game between the aggregator and consumers is a subset of the set of possible outcomes of a Nash Bargaining Game, where consumers bargain for fair benefit allocations. Nash-Cournot games can model the power market competition where the bidding prices depend on the bidding quantity \cite{li_risk-averse_2017, li_distributed_2023}. However, the linear bidding curve could be naive for ISO-type auctions \cite{hobbs_strategic_2000-1}, and there are no overbidding regulations to secure system reliability in Nash-Cournot games.}

The risks in power markets exist in the energy price volatility \cite{li_distributed_2023}, changing risk appetites \cite{vespermann_risk_2021}, and uncertain renewable energy outputs  \cite{jiang_flexibility_2022}. This will impact the strategy set of players and the properties of game equilibrium \cite{vespermann_risk_2021}. In the literature, the player's strategies are modeled using risk-aware formulations, including the chance-constrained (CC) \cite{bruninx_interaction_2020}, conditional value-at-risk (CVaR) \cite{ordoudis_energy_2021}, and distributionally robust formulations \cite{ding_distributionally_2022}. These formulations assume an underlying distribution or a family of distributions for uncertain variables and model the expected market outcomes at a certain chance. Since a game result is the outcome of all players regardless of game types, a joint chance-constrained (JCC) formulation is necessary to model the occurrence of uncertain market outcomes due to multiple players' strategies jointly \cite{li_data-driven_2022, mazadi_impact_2013}. 

The JCC problem is generally intractable. Their solutions can be approximated by decomposing the JCC into tractable single CCs with predetermined violation rates by enforcing \textit{Boole's inequality}. This means that the chance that at least one of the events happens is no greater than the sum of chances of all individual events \cite{xie_optimized_2022}. A particular case is the \textit{Bonferroni approximation}, which assumes that these CCs are entirely independent and the sum of individual violation rates is equal to the joint violation rate  \cite{ding_distributionally_2022, xie_optimized_2022}. However, CCs with the same uncertainty source (e.g., renewable power sources within an area) are usually positively correlated. This assumption neglects their correlation and increases solution conservativeness. Ref. \cite{xie_optimized_2022} proves the solution conservativeness of the JCC problem will increase with the number of decomposed CCs based on the \textit{Bonferroni approximation}. Ref. \cite{chen_cvar_2010} proves that if these constraints are entirely correlated, the efficacy of the \textit{Bonferroni approximation} will diminish. 

Although the correlation of decomposed CCs impacts solution conservativeness in a JCC problem, this topic has received very limited discussions in game designs. Ref. \cite{mazadi_impact_2013} first used the JCC formulation in a Nash-Cournot game to simulate the impact of uncertain wind power outputs in the electricity market considering the congestion arbitrage. However, it simply assumes all the wind farms are independent and have the same confidence levels in predictions. Recent work \cite{li_data-driven_2022} uses a DRJCC formulation to model the generalized Nash equilibrium problem for coordinating the renewable energy aggregators in the local electricity market. The work considers the maximum power outputs of renewable energy aggregators as an uncertain distribution. The model solves their bidding strategies and market equilibrium using the best response algorithm. Nevertheless, the paper does not discuss the correlation between players, considering the game is simulated with varying numbers of players. 

To identify the correlation between CCs and reduce solution conservativeness, a common method is pre-solving the JCC with a conservative guess and progressively tightening the worst-case bound. Ref. \cite{chen_cvar_2010} proposed a tractable CVaR formulation to model the worst-case bound of the JCC problem. An auxiliary variable is introduced to tighten the bound. This optimized CVaR method has been applied in the energy and reserve joint dispatch problem \cite{ordoudis_energy_2021}. Refs. \cite{baker_joint_2019} and \cite{jia_iterative_2021} solve the JCC optimal power flow (OPF) problem. The overlapped violation events due to the correlated CCs are estimated analytically and eliminated iteratively. Recent work \cite{arrigo_embedding_2022} incorporates the correlation in wind power generations into the DRO as a constraint. All these works \cite{baker_joint_2019, jia_iterative_2021, arrigo_embedding_2022} effectively reduce solution conservativeness but require a comprehensive analysis of violation events or the specified correlation, which is hard to obtain in a changing power market. Recent works \cite{inatsu2022bayesian, shapiro_bayesian_2023} employ the Bayesian method to estimate unknown parameters of the underlying distributions for DRO and DRCC problems. However, no attempts are made to solve the JCC or DRJCC problems using Bayesian optimization considering the unspecific correlation.

This paper designs a DRJCC game framework to coordinate solar-powered microgrids in the energy and reserve markets. Compared with our early work in ref. \cite{ding_coordinating_2023} \footnote {Some early results of this work have appeared in the paper of the 14th ACM International Conference on Future Energy Systems \cite{ding_coordinating_2023}}, which uses the DRCC formulation to model the chance of breaching individual contracts, this work considers the joint contract violation rate across the power network for system reliability. We also develop a novel end-to-end Bayesian approach to approximate the optimal solution of the proposed DRJCC game model, considering the charging correlation among players. Our contributions are two-fold,

(i) From the perspective of market coordination, we consider the risk-aware strategies of both the system operator and microgrids and their joint impacts in a market game. A two-stage DRO formulation models the risk-aware bidding of microgrids, and the ambiguity set is constructed by empirical solar forecast errors in the California power market. A JCC-based system regulation is imposed on the system operator’s side to manage the reserve performance in the real-time market. The system operator’s and multiple microgrids’ interactions are simulated in a Stackelberg leader-follower game.

(ii) From the perspective of algorithm design, the proposed DRJCC model becomes intractable if the individual violation rate of decomposed DRCCs is a design variable \cite{nemirovski_convex_2007}. On the condition of the existence and uniqueness of the game solution, a novel end-to-end Bayesian optimization method is proposed to approximate individual contract violation rates and solve the game equilibrium. Compared with the latest works about DRCC problems \cite{inatsu2022bayesian,  shapiro_bayesian_2023}, we first attempt to use Bayesian optimization for tackling the DRJCC problem. Results show that the optimal individual violation rate can avoid overbidding or underbidding while securing reliable reserve services.

The rest of this paper is organized as follows. Section III presents the proposed DRJCC game types and players. Section IV demonstrates the model formulation based on MPEC. Then, the Bayesian optimization algorithm for solving the DRJCC problem is introduced in Section V. Sections VI and VII present case studies and numerical results. Section VIII concludes the paper. 

\section{DRJCC market game considering uncertain reserve provision}

We model the two-settlement electricity market, which consists of a day-ahead forward market and a real-time market. We use a bilevel optimization framework to model the interactions of ISO (upper level) and renewable microgrids (lower level). Rather than an aggregator as a trading intermediary between the upstream wholesale market and downstream service providers, the role of each microgrid is similar to a DSO in this paper. They schedule the power dispatch of individual microgrids and bid in the upstream markets.

\begin{figure}[!h]
    \centering
    \includegraphics[width=3in]{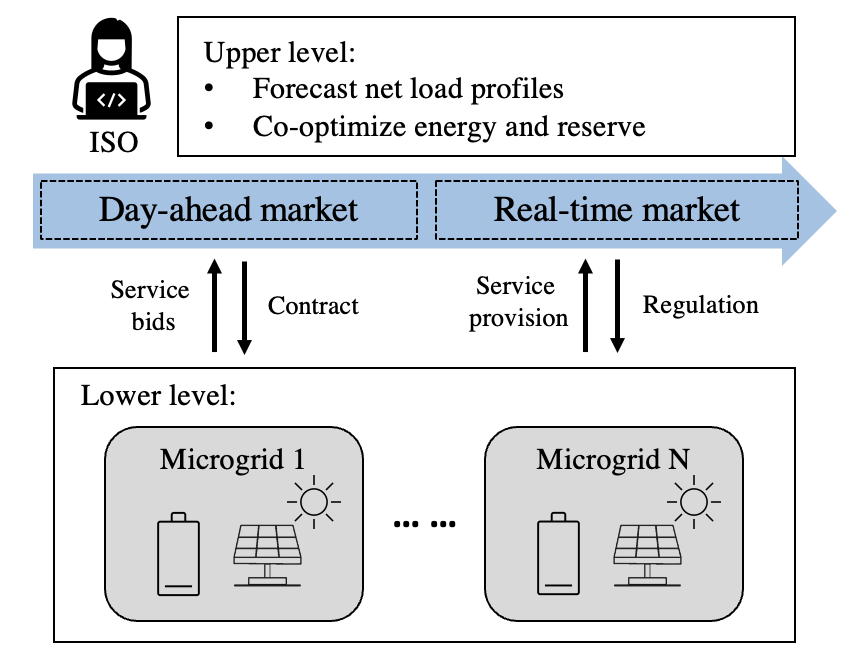}
        \caption{\textcolor{black}{Bilevel Stackelberg game framework for coordinating renewable microgrids in two-settlement electricity markets}}
    \label{schematic}
\end{figure}

This game is a one-leader, multiple-follower Stackelberg game and includes two types of players. ISO is a leader, and contracted microgrids are followers. As shown in Fig. \ref{schematic}, in the day-ahead market, ISO sets energy and reserve requirements using probabilistic net load forecasts (i.e., power demands subtracting solar power generations) at each network node. Generators and renewable microgrids submit bids and sign day-ahead services contracts. ISO knows all bids and conducts market clearing. Each microgrid does not learn other players’ bidding strategies and is assumed not to affect the market.

However, due to the uncertainty of renewable generation, reserve services from renewable microgrids might under-deliver in real time. To guarantee system reliability, ISO will assess the reserve provision performance, regulate contract violations in the market, and broadcast the maximum individual contract violation rates to microgrids. Meanwhile, each microgrid can leverage the historical forecast data and adjust its service bidding to avoid penalties or payment reductions due to under-delivered services.

%\subsection{Game type}
%
%This game is a one-leader-multiple-follower Stackelberg game. In this game, the leader anticipates followers’ actions and moves first, while followers take moves afterward only based on the leader’s actions. In this case, ISO knows all the microgrid bids (e.g., capacity, penalty, and risk appetite) and conducts the market clearing. Each microgrid bids based on the market clearing result without knowing other players’ bidding strategies. They are assumed to be in perfect competition and not to affect the wholesale market.

\section{Model formulation}

We first formulate the risk-aware optimization problem of ISO and microgrids. To model their interactions, we integrate these optimization models into a complementarity model that can represent the simultaneous optimization problems of one or several decision-makers \cite{ruiz_tutorial_2014}.

\subsection{Wasserstein-metric ambiguity set to model the worst-case reserve performance}

%Since microgrids are mainly supplied by renewable power sources such as solar and wind, their real-time reserve provision performances are modeled based on historical renewable forecast error samples.

Since microgrids are mainly supplied by solar power, we use historical renewable forecast error samples to construct a Wasserstein-metric ambiguity set and model the worst-case reserve performance. It is a moment-free method without any underlying distribution assumptions leveraging the empirical data. A ball-shape set space is constructed, which allows adjusting conservativeness by its radius. We chose the Wasserstein distance for its advantages, including symmetry and capturing the worst-case distribution. \cite{gao_distributionally_2023}. \footnote{\textcolor{black}{Other moment-free distances, such as KL-divergence, can also be employed theoretically. However, the KL-divergence is an asymmetric metric, meaning that the divergence from the distribution $\mathbb P$ to $\mathbb Q$ is not equal to the divergence from the distribution $\mathbb P$ to $\mathbb Q$ \cite{lovric_kullback-leibler_2011}. This will bring more complexity in defining conservativeness.}}

Let the error sample vector $ \boldsymbol s := \{s_1, s_2,..., s_n\}$ be a vector including $n$ random samples of the uncertain variable $ \zeta$ at time $t$. If $\boldsymbol s_1 \sim \mathbb P_1 $ and $\boldsymbol s_2 \sim \mathbb P_2 $, the Wasserstein metric $W (\mathbb P_1, \mathbb P_2) $ between two distributions $\mathbb P_1$ and $\mathbb P_2$ is given by \cite{netessine_wasserstein_2019},

\begin{equation}
\label{ambiguity set_2}
W (\mathbb P_1, \mathbb P_2) = \inf_{ \mathcal Q \in \mathcal R (\Xi)} \{  \int ||s_1 - s_2 || \mathcal Q ( d \boldsymbol s_1,  d  \boldsymbol s_2) \}
\end{equation}

The term $||s_1 - s_2 ||$ is the distance between two samples, and $\mathcal Q ( d \boldsymbol s_1,  d  \boldsymbol s_2)$ is the joint distribution with marginal distributions $\mathbb P_1$ and $\mathbb P_2$ on the support $\mathcal R (\Xi)$. The Wasserstein-metric ambiguity set $\mathcal P$ is a ball space of radius $\nu \geq 0$ concerning Wasserstein distance, centered at a prescribed reference distribution $\mathbb P_{\zeta}$ based on samples $\boldsymbol s$. 

\begin{equation}
\label{ambiguity set_2}
\mathcal P := \{\mathbb P : \mathbb P \in \mathcal R (\Xi),W (\mathbb P, \mathbb P_{\zeta}) \leq \nu \}
\end{equation}

\begin{figure}[!h]
    \centering
    \includegraphics[width=3.5in]{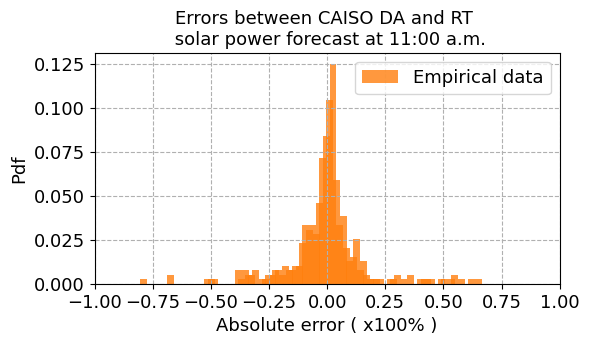}
        \caption{The solar power forecast error distribution between day-ahead and real-time markets at 11:00 a.m. \cite{california_iso_oasis_nodate}}
    \label{error}
\end{figure}

This work uses historical error samples of solar power forecasts from California Independent System Operator (CAISO) during 2021 and 2022 \cite{california_iso_oasis_nodate}. We choose two real-time market intervals at 15 and 45 mins in each hour to calculate these errors between the day-ahead and real-time markets. Fig. \ref{error} shows the forecast error distribution at 11:00. For each time interval, we choose $N_S = 200$ to construct the ambiguity set. These error samples are integrated into the model as a matrix $\textbf s \in \mathbb R^{N_S \times N_T}$.

\subsection{ISO's optimization for the energy and reserve joint dispatch}
  
This work considers a power network with $N_B$ buses and $N_{LN}$ lines. Excluding the slack bus, the bus susceptance matrix denotes as $\textbf{B}_{bus} \in \mathbb R^{(N_B-1)  \times (N_B-1)}$, and the flow susceptance matrix denotes as $\textbf{B}_{line} \in \mathbb R^{N_{LN} \times (N_B-1)}$. The power system incorporates a set of thermal power plants $k\in \mathcal K$, net loads $j \in \mathcal J$, and microgrids $i \in \mathcal I$. The positions of $N_G$ thermal power plants, $N_L$ net loads, and $N_{MG}$ microgrids are mapped by sparse matrices $\textbf{C}^{g} \in \mathbb R^{N_B \times N_G }$, $\textbf{C}^l \in \mathbb R^{N_B \times N_L}$ and $\textbf{C}^{mg} \in \mathbb R^{N_B \times N_{MG}}$. \textbf{Bold letters} represent the vector of a variable or a parameter across the network over the optimization horizon (e.g., $\textbf{P}_{nl}$ is all net loads over the period $\mathcal T$). If the real-time market clearing interval is $\Delta t \in \mathcal T$, the multi-period power dispatch is as follows.

1) \textit{The objective function:} ISO's cost function $ \mathcal  J^{H}$ consists of two parts. The first part is the cost of energy and reserve provision of thermal power plants, assuming they have the quadratic function in the economic dispatch. The second part is the net procurement costs for energy and reserve services from microgrids. This is calculated as the product of bidding price $\pi_i^r, \pi_i^e$ and bidding quantity $R^{up}_{t, i}, P^{ex}_{t, i}$.

\begin{equation} 
\begin{aligned}
(M1)  \min_{\Xi_{hl}} \mathcal  J^{H} =  \min \{  \sum_{t\in \mathcal T} \sum_{k\in \mathcal K}  (m_k^g( \Delta t P^g_{t, k})^2 +  \Delta t m_k^{gr}R^g_{t, k})  \\
+   \sum_{t\in \mathcal T} \sum_{i\in \mathcal I}  \Delta t (\pi_i^r R^{up}_{t, i}  -  \pi_i^e P^{ex}_{t, i}) \}
 \label{objective1}
 \end{aligned}
 \end{equation}
 
\textcolor{black}{ where the variable set is $\Xi_{hl} = \{P^g,  R^g, \pi^r, R^{up}, \pi^e, P^{ex} \}$}
 
2) \textit{Power and reserve provision}: 

\begin{equation} 
 \sum_{k \in \mathcal K}  P^{g}_{t, k} + \sum_{i \in \mathcal I}  P^{ex}_{t, i} = \sum_{j \in \mathcal J}{P^{nl}_{t, j} \quad t \in \mathcal T}\label{power_balance}\\
 \end{equation}
 
   \begin{equation} 
\sum_{k \in \mathcal K} R^g_{t, k}  + \sum_{i \in \mathcal I} R^{up}_{t, i} \geq R^{req}_t \quad t \in \mathcal T
  \label{reserve_balance}
 \end{equation}
 
   \begin{equation} 
  \textbf P^g + \textbf R^g \leq \overline {\textbf P^g} \label{reserve_upper}\\
 \end{equation}
 
  \begin{equation} 
 \textbf P^g \geq 0, \textbf R^g \geq 0 \label{non-negative}\\
 \end{equation}
 
Constraints (\ref{power_balance}) and (\ref{reserve_balance}) are for power balance and reserve adequacy. Constraint (\ref{reserve_upper}) regulates the sum of energy and reserve provisions is less than the power capacity of thermal power plants, and service provisions are non-negative (\ref{non-negative});

  \begin{equation} 
\begin{aligned}
\textbf{P}^{inj} =\textbf{C}^{g} (\textbf{P}^{g}+ \textbf R^{g}) - \textbf{C}^{mg} (\textbf{P}^{ex} - \textbf{R}^{up})- \textbf{C}^{l} \textbf{P}^{nl} \label{power_injection}
\end{aligned}
 \end{equation}

 \begin{equation}
|| \textbf{B}_{line}\textbf{B}_{bus}^{-1} \textbf{P}^{inj} || \leq \overline {\textbf{L}} 
\label{line_limit_1}
 \end{equation}
 
 The net power injection is denoted as (\ref{power_injection}), and the power flow should be within the line power capacity (\ref{line_limit_1}).

3) \textit{Reserve performance regulations based on service contracts}: In the real-time market, ISO assesses reserve performance and regulates the occurrence of contract violations for system reliability. We formulate this regulation as a JCC based on the Wasserstein-metric ambiguity set $\mathcal P$. Assuming the maximum undelivered reserve capacity is $\overline {R^{fup}}$ in each day-ahead contract, the JCC regulates the possible violations of all contracts less than the chance of $\epsilon^{jot}$.

  \begin{equation} 
\inf_{\mathbb P \in \mathcal P }\mathbb P ( \bigcap_{i \in \mathcal I} \{\sum_{t\in\mathcal T} \Delta t \zeta R_{i, t}^{up} \leq \overline {R^{fup}}\}) \geq 1-\epsilon^{jot}
\label{DRJCC}
 \end{equation}
 
Where $\epsilon^{jot}$ is predetermined and depends on the availability of the system contingency reserve. To assess the performance of individual microgrids, the JCC (\ref{DRJCC}) is decomposed into single CCs.

  \begin{equation} 
\inf_{\mathbb P \in \mathcal P }\mathbb P ( \sum_{t\in\mathcal T} \Delta t \zeta R_{i, t}^{up} \leq \overline {R^{fup}}) \geq 1-\epsilon^{ind}_i \quad i \in \mathcal I
\label{DRCC}
 \end{equation}
 
Where $\epsilon_i^{ind}$ represents the maximum chance of violating each microgrid's contract. As the correlation among players is unknown, $\epsilon_i^{ind}$ could not be solved directly. Instead, we could estimate its value within a range.
 
 \begin{equation} 
\epsilon_i^{ind} \in [\frac{\epsilon^{jot}}{N_{MG}},  \epsilon^{jot}]
\label{range}
 \end{equation}
 
\textbf{Remark 1:} \textit{The joint and individual violation rate needs to satisfy Bonferroni’s inequality $\sum^{N_{MG}}_{i=0}\epsilon^{ind}_i \leq \epsilon^{jot}$. In this case, the most conservative choice of  $\epsilon^{ind}_i$ is $\epsilon^{ind}_i  =\frac{ \epsilon^{jot}}{N_{MG}}$ for all $i = 1, . . . , N_{MG}$ \cite{nemirovski_convex_2007}, which is set as the lower bound and also known as \textit{Bonferroni approximation}. This case assumes no correlation among microgrids (i.e., decomposed CCs). The least conservative choice is $\epsilon^{ind}_i$ is $\epsilon^{ind}_i  = \epsilon^{jot}$ which is set as the upper bound. This assumes all microgrids are completely correlated and have the same performance. }

5) \textit{CVaR approximation and reformulation:} As the proof (Theorem 2.2) in ref. \cite{zymler_distributionally_2013}, each DRCC (\ref{DRCC}) can be approximated as a CVaR constraint and then reformulated into the following convex constraints based on refs. \cite{mohajerin_esfahani_data-driven_2018, arrigo_wasserstein_2022}.

  \begin{equation} 
 -\eta^u_i - \frac{1}{\epsilon_i^{ind}} ( \nu l^u_{i} + \frac{\boldsymbol 1^{\top}\textbf W^u_i}{N_S} ) \geq 0
\label{convex_constraint_1}
 \end{equation}
 
   \begin{equation} 
 \textbf W^u_i -    \Delta t \textbf s (\textbf R^{up}_ i)^{\top}+ \overline {R^{fup}} + \eta^u_i \geq 0
\label{convex_constraint_2}
 \end{equation}
 
     \begin{equation} 
 \textbf W^u_i \geq 0
\label{convex_constraint_3}
 \end{equation}
 
     \begin{equation} 
|| \Delta t \textbf s (\textbf R^{up}_ i)^{\top}|| \leq  l^u_i
\label{convex_constraint_4}
 \end{equation}
 
Where $\textbf R^{up}_i \in \mathbb R^{1 \times N_T}$ is the reserve vector of microgrid $i$, and variables $\eta^u$, $ l^u$, $ \textbf W^u \in \mathbb R^{N_S \times 1}$ are auxiliary variables.

  \subsection{Data-driven bidding strategy considering the worst-case reserve performance}
  
Microgrids conduct risk-aware bidding in a two-settlement power market and are modeled as a two-stage DRO formulation. The first stage models the energy and reserve bidding in the day-ahead market, and the second stage models the expected penalty of under-delivered reserve considering the worst-case reserve performance. This penalty is modeled as involuntary load shedding depending on the first stage's empirical error samples and day-ahead reserve bids. The following subsection presents the model formulation, and the dual variable of each constraint is presented in the left bracket. 
  
 1) \textit{The objective function:} Each microgrid $i\in \mathcal I $ tries to minimize its cost $\mathcal  J^{L}$, which consists of three parts. The first part is the energy and reserve provisions cost, and the second is the revenue from energy and reserve service bidding. The final part models the monetary penalty due to the unavailable reserve service. We assume the penalty is the worst-case load-shedding loss and use the DRO formulation based on the Wasserstein-metric ambiguity set $\mathcal P$. 

\begin{equation} 
\begin{aligned}
(M2) \min_{\Xi_{ll}} \mathcal  J^{L} = \min  \{ \sum_{t\in\mathcal T} (m^s_i ( \Delta t P^s_{i, t}) ^2 + \Delta t m^{mgr}_i R^{up}_{i, t} \\
 +\Delta t \pi^e_iP^{ex}_{i, t} - \Delta t \pi^r_i R^{up}_{i, t} )\} \\
 + \sup \mathbb E_{\mathbb P \in \mathcal P} [  \Delta t  \zeta R^{up}_{i, t} V_i] 
 \label{objective2}
 \end{aligned}
 \end{equation}
 
where the variable set is $\Xi_{ll} = \{ P^s, \pi^r, R^{up}, \pi^e, P^{ex}\}$. The load-shedding loss due to the under-delivered reserve is calculated by multiplying each microgrid's under-delivered reserve capacity and the VoLL ${V_i}$. 
 
 2) \textit{Convex reformulation of DRO:} The DRO term in the objective function (\ref{objective2}) is reformulated based on Corollary 5.4 in ref. \cite{mohajerin_esfahani_data-driven_2018}.
 
    \begin{equation} 
 \sup  \mathbb E_{\mathbb P \in \mathcal P} [ \sum_{t\in\mathcal T} \Delta t\zeta R^{up}_{i,t} V_i] =  \nu l^l_i + \frac{\boldsymbol 1^{\top}\textbf W^l_i}{N_S}
\label{lower_0}
 \end{equation}

  \begin{equation} 
(\boldsymbol \mu^{dro1})  \quad \textbf W^l_i  - \Delta t \textbf s (\textbf R^{up}_i)^{\top} V_i + \eta^l_i \geq 0 \label{lower_1}
 \end{equation}
 
   \begin{equation} 
(\boldsymbol \mu^{dro2})  \quad \textbf W^l_i \geq 0 \label{lower_2}
 \end{equation}
 
    \begin{equation} 
(\boldsymbol \mu^{dro3})  \quad \Delta t\textbf s (\textbf R^{up}_i)^{\top} V_i + l^l_i \geq 0 \label{lower_3}
 \end{equation}
 
    \begin{equation} 
(\boldsymbol \mu^{dro4})  \quad \Delta t\textbf s (\textbf R^{up}_i)^{\top} V_i -  l^l_i   \leq 0 \label{lower_4}
 \end{equation}
 
 where $ l^l $, $\eta^l$ and $ \textbf W^l \in \mathbb R^{N_S \times 1}$ are auxiliary variables in the reformulation. The dual variables $\boldsymbol \mu^{dro1}$,$ \boldsymbol \mu^{dro2}$, $\boldsymbol \mu^{dro3}$ and $\boldsymbol \mu^{dro4}$ all have the same dimension $\mathbb R^{N_S \times 1}$.
   
 3) \textit{Power balance considering exchange with the upstream grid}: Since the dispatch in the microgrid is unknown for ISO during the market-clearing, microgrids are modeled as flexible power sources $P^s_{i, t}$, which is a proxy of renewable generation and energy storage.

\begin{equation} 
(\mu^{pb})  \quad  P^{ex}_{i, t} + P^s_{i, t} =  P^{nl}_{i, t} \quad t \in \mathcal T \label{power_balance_LL}
 \end{equation}
 
4) \textit{Energy and reserve bids}: Each microgrid submits energy and reserve bids in the day-ahead forward market. The capacity range of the energy service bidding is given by,
 
 \begin{equation} 
(\mu^{dnpl}, \mu^{uppl})  \quad - \underline{P^s} \leq P^s_{i,t} \leq \overline{P^s} \quad t \in \mathcal T
\label{solar_energy_bid}
 \end{equation}
 
The capacity range for the reserve service bidding is between zero and a fraction of the flexible power generation level $ \gamma$. 

 \begin{equation} 
( \mu^{dnrl},  \mu^{uprl}) \quad 0 \leq  R^{up}_{i, t} \leq \gamma P^{s}_{i, t} \quad t \in \mathcal T
\label{solar_reserve_bid}
 \end{equation}

\subsection{Complementarity model for market equilibrium}

lSO and microgrids' optimization are interlinked by energy and reserve bids (i.e., $\pi^r R^{up}  - \pi^e P^{ex}$). To form a one-leader-multiple-follower problem, the MPEC technique integrates two optimization problems to obtain market equilibria. (The best response method cannot be used here as it might get multiple near-optimal points, demonstrated in \textbf{Remark 2}). 
Based on the strong duality and KKT rule, the dual problem of model M2 for $N_{MG}$ of microgrids is integrated into the upper-level model M1. We omit the indexes $i, t$ for conciseness. $\textbf1^s$ represents a unit vector with the same dimension as the error sample vector $\boldsymbol s$. First, stationary constraints are formulated for six variables $P^{s}, P^{ex},R^{up}, l^l, W^l$ and  $\eta^l$ in model M2.

   \begin{equation} 
2  \Delta t m^s P^s + \mu^{pb}  -  \mu^{dnpl} + \mu^{uppl} -  \gamma \mu^{uprl} = 0    \label{ps}
 \end{equation} 
 
  \begin{equation} 
\boldsymbol \pi^e +  \boldsymbol \mu^{pb}  = 0    \label{pex}
 \end{equation} 
 
   \begin{equation} 
   \begin{aligned}
  \Delta t ( m^{mgr} -   \pi^r) + \Delta t \boldsymbol s^{\top} ( \boldsymbol \mu^{dro1} -  \boldsymbol \mu^{dro3} +  \boldsymbol \mu^{dro4})  {V} \\ - \mu^{dnrl} +  \mu^{uprl} = 0   \label{rup}
 \end{aligned}
 \end{equation}  
 
    \begin{equation} 
 \nu - (\textbf1^s)^{\top} \boldsymbol \mu^{dro3} -   (\textbf1^s)^{\top} \boldsymbol \mu^{dro4} = 0   \label{pex}
 \end{equation}  
 
  \begin{equation} 
\frac{1}{N_S} - \boldsymbol \mu^{dro1} - \boldsymbol \mu^{dro2} = 0 \label{W}
 \end{equation} 
 
   \begin{equation} 
 (\textbf1^s)^{\top}  \boldsymbol \mu^{dro1}  = 0 \label{eta}
 \end{equation} 
 
Complementary slackness constraints are,
 
  \begin{equation} 
  0 \leq  (\textbf W^l  -  \Delta t \textbf s (\textbf R^{up})^{\top} V + \eta^l) \bot \boldsymbol \mu^{dro1} \geq 0  \label{csc1}
  \end{equation} 
  
    \begin{equation} 
  0 \leq \textbf W^l   \bot \boldsymbol \mu^{dro2} \geq 0  \label{csc2}
  \end{equation} 
  
    \begin{equation} 
  0 \leq  ( \Delta t \textbf s (\textbf R^{up})^{\top} V  + l^l ) \bot \boldsymbol \mu^{dro3} \geq 0  \label{csc3}
  \end{equation} 
  
    \begin{equation} 
  0 \leq  (l^l -  \Delta t \textbf s (\textbf R^{up})^{\top} {V})  \bot \boldsymbol \mu^{dro4} \geq 0  \label{csc4}
  \end{equation} 
  
   \begin{equation} 
0 \leq  (\overline{ P^s} -  P^s) \bot \mu^{uppl} \geq 0   \label{csc5}
 \end{equation} 
  
      \begin{equation} 
  0 \leq (\underline{ P^s} + P^s)  \bot  \mu^{dnpl} \geq 0  \label{csc6}
  \end{equation}

 \begin{equation} 
0\leq R^{up} \bot \mu^{dnpl} \geq 0   \label{csc7}
 \end{equation} 
 
   \begin{equation} 
0 \leq (\gamma P^{s} - R^{up} ) \bot \mu^{uppl} \geq 0  \label{csc8}
 \end{equation}

Constraints (\ref{csc1}) - (\ref{csc8}) are non-linear. We thus use the big-M technique. Assuming $a$ and $b$ are variables, the non-linear constraint on the left-hand side can be transformed into the four convex constraints on the right-hand side.

   \begin{equation} 
0\leq a \bot b \geq 0 \Rightarrow a \geq 0, b\geq 0, a \leq MU, b \leq M(1-U)   \label{bigm}
 \end{equation} 
 
 where $M$ is a number far larger than $a$ and $b$, and $U$ has the same dimension as $a$ and $b$.
 
The objective function of ISO's optimization problem is nonlinear after MPEC reformulation, with the multiplication of price and capacity variables. We linearize it based on the strong duality theory; the details can be found in ref. \cite{ding_coordinating_2023}. The reformulated model is given by,

\begin{equation} 
\begin{aligned}
(M3) \quad \min_{\Xi_{bi}} \mathcal  J = \min\{ \sum_{t\in \mathcal T}\sum_{k\in \mathcal K}  ( m_{k}^g ( \Delta t P^g_{t,k})^2 +  \Delta t m_{t, k}^{gr}R^g_{t, k}) \\
+ \sum_{t\in \mathcal T} \sum_{i\in \mathcal I} [m_i^s (\Delta t P_{t, i}^s) ^2 +  \Delta t  m_i^{mgr} R_{t, i}^{up}  + (\nu l^l_i + \frac{\boldsymbol 1^{\top}\textbf W^l_i}{N_S}) \\
- \Delta t \mu^{pb}_{t, i} {P_{t, i}^l} + \Delta t \mu^{uppl}_{t, i} \overline{P^s} +  \Delta t \mu^{dnpl}_{t, i} {\underline{ P^s}} ]\}\\
\qquad \textrm{s.t. \quad equations} \quad (\ref{power_balance}) - (\ref{convex_constraint_4}) \quad \textrm{and} \quad (\ref{ps}) - (\ref{csc8})
  \label{objective3}
 \end{aligned}
 \end{equation}

 where the variable set ${\Xi_{bi}}$ includes $P^g$, $R^g$, $P^s$, $R^{up}$, $R^{ex}$, $\pi^r$, $\pi^e$, $l^l $, $\eta^l$, $ \textbf W^l $, $\eta^u$, $ l^u$, $ \textbf W^u$, $\epsilon^{ind}$, $ \mu^{dro1}$,  $\mu^{dro2}$,  $\mu^{dro3}$, $\mu^{dro4}$, $\mu^{pb}$, $\mu^{dnpl}$, $\mu^{uppl}$,
 $\mu^{dnrl}$, $\mu^{uprl}$. 
 
\section{Approximating the optimal solutions using Bayesian optimization}
   
The reformulated MPEC model M3 is intractable with unknown individual violation rates $\boldsymbol  \epsilon^{ind}$. We thus introduce Bayesian optimization to approximate the solution iteratively, a methodology designed to approximate black-box functions. Its applications include hyper-parameter tuning for deep learning models \cite{masaki_adachi_sober_nodate} and solving combinatorial optimizations \cite{shahriari_taking_2016}. It best suits the expensive-to-evaluate functions and tolerates the stochastic noise in function evaluations \cite{shahriari_taking_2016}. Here, the target black-box function is the unspecified relationship between the individual violation rates $\boldsymbol \epsilon^{ind}$ and the empirical joint violation rate $\epsilon^{e}$, denoted as $\epsilon^{e} = \mathcal F (\boldsymbol \epsilon^{ind})$. Its evaluation requires out-of-sample tests and is considered noisy observations.

\textbf{Remark 2:} \textit{The prerequisite of using global optimization to solve a game is the existence and uniqueness of the game solution. The stochastic game model M3 must have unique equilibria to ensure the black-box function is unchanged during the solution exploration. We prove it in Appendix A. }

%\textcolor{black}{Bayesian optimization is thus introduced to solve the black-box function $\mathcal F $, given the following two reasons. First, the black-box function cannot be explicitly computed. Its evaluation process requires out-of-sample tests, which are considered to have noisy observations. Second, during the optimization process, we do not consider data transfer costs (e.g., ISO broadcasts violation rate and prices), which can be substantial if the power network scales up. Bayesian optimization can guarantee performance within a limited number of iterations and reduce evaluation costs. The following subsections summarize its algorithm, the design of the evaluation, and acquisition functions in the process. }

\subsection{Summary of the algorithm}

\begin{algorithm}
\SetAlgoLined
Initialize model M3 and define the exploration space; \\
Evaluate the initial samples $\boldsymbol \epsilon_o$; \\
\For{$ m = 1,....., N_{IT}$  }{
      Update the surrogate model with observations $ \textbf O_{m}$\;
      Identify new samples $\boldsymbol \epsilon_{m+1} = \underset{\epsilon}{\mathrm{argmax}}  \mathcal G (\boldsymbol \epsilon_m;  \textbf O_m)$\;  
      Query the evaluation function $ \mathcal H (\boldsymbol \epsilon_{m+1})$\;
     \textcolor{black}{ Augment observations $ \textbf O_{m+1} := \{ \textbf O_{m}; (\boldsymbol \epsilon_{m+1}, \epsilon^e_{m+1})\}$\;}
    }
\caption{Solving M3 using Bayesian optimization}
\end{algorithm}

As in Algorithm 1, we initialize the problem to solve, model M3, and define the exploration space of the individual violation rates (\ref{range}). The dimension of the sample vector is $N_{MG}$, and the bounds ensure that the black-box function is feasible and continuous over the space. Then, a Gaussian process regressor (GPR) is used as the surrogate model for fitting the unknown block-box function $\mathcal F$ within $N_{IT}$ iterations.  \textcolor{black}{GPR assumes a prior distribution for the unobserved value $ \boldsymbol \epsilon^e_m $. Given the distributions of observations $\mathbb P ( \textbf O_m)$ and the likelihood $\mathbb P ( \textbf O_m| \boldsymbol \epsilon^e_m ) $ (i.e., the maximum difference between the estimation and the ground truth), the posterior of GPR (i.e., the probability distribution of the approximated value) is, }

  \begin{equation} 
\mathbb P (\boldsymbol \epsilon^e_m | \textbf O_m) = \frac {\mathbb P ( \textbf O_m| \boldsymbol \epsilon^e_m ) \mathbb  P ( \boldsymbol \epsilon^e_m )}{\mathbb  P ( \textbf O_m)}
\label{gpr}
 \end{equation}

In each iteration $m$, the posterior is updated using all observations $\textbf O_{m}$. After that, new samples are proposed for the next iteration using the acquisition function $ \mathcal G (\boldsymbol \epsilon_m; \textbf O_m)$. These samples are evaluated in the function $\mathcal H(\boldsymbol\epsilon_{m+1})$ and added to the observation $\textbf O_{m+1} $. This procedure is executed for $N_{IT}$ iterations. 

\subsection{Evaluation of the black-box function}

Bayesian optimization finds the optimal individual violation rates $\boldsymbol \epsilon^{ind}$ so that the empirical outcome $\epsilon_m^{e}$ is as close as to the requirement $\epsilon^{jot}$. As the black box function $\mathcal F$ cannot be explicitly computed, an evaluation function $\mathcal H(\boldsymbol \epsilon_m)$ assesses results in each iteration. This is the absolute difference between the empirical and predefined joint violation rates.

  \begin{equation} 
\mathcal H(\boldsymbol \epsilon_m) = |\epsilon_m^{e} - \epsilon^{jot}| 
\label{evaluated_function}
 \end{equation}

Where the empirical joint violation rate $\epsilon_m^{e}$ is obtained from $N_{OS} = 150$ independent out-of-sample tests, In each test, we use an error sample $s$ different from those in the ambiguity set and calculate the under-delivered reserve $ R^{under}_s:= \Delta t s (\textbf R^{up})^{\top}$  for each microgrid. A violation event is when any microgrid violates the reserve regulation (i.e., constraint (\ref{DRCC})). The empirical joint violation rate is calculated as the ratio of violation events to the total number of tests.

 \begin{equation} 
\epsilon_m^{e}:= \frac{\sum_{s=0}^{N_{OS}}\mathbbm{1}_{(R^{under}_s > \overline{R^{fup}})}}{N_{OS}}
\label{empirical_violation}
 \end{equation}
 
 \subsection{Sample selection using the acquisition function}
 
The acquisition function decides the new sample in the iteration $m$. Different acquisition functions include the lower confidence bound, the probability of improvement, and the expected improvement \cite{wilson_maximizing_2018}. We choose the expected improvement after experiments considering the explore-exploit tradeoff. The function is the maximum expected improvement at the sample $\boldsymbol \epsilon_m$ estimated from all out-of-sample tests. 
 
    \begin{equation} 
 \mathcal G (\boldsymbol\epsilon_m;  \textbf O_m) = \mathbb E [\mathcal H^{'}-\mathcal H(\boldsymbol \epsilon_m)] ^{+}
\label{out_of_sample}
 \end{equation}
 
Where $[ \cdot ]^{+} = \max\{ \cdot, 0 \}$, and $\mathcal H^{'}$ denotes the minimum value obtained in the evaluation function. Equation (\ref{out_of_sample}) defines the maximum expected improvement based on the posterior distribution in each iteration, which is greater than zero. Then, a new sample in the next iteration is selected, which gives the maximum value of (\ref{out_of_sample}).

   \begin{equation} 
\boldsymbol \epsilon_{m+1} = \underset{\epsilon}{\mathrm{argmax}} \mathcal G (\boldsymbol \epsilon_m ; \textbf O_m)
\label{new_sample_point}
 \end{equation}
 
\section{Case studies} 

Case studies are developed using parameters in Table \ref{System_para_2} and IEEE 30-bus network \cite{power_systems_test_case_archive_30_nodate}. The power network consists of passive loads and microgrids at 24 load buses. Each microgrid has a 2MW load on average, a 2MW solar PV generation, and a 3MW/6MWh ES system. 

\begin{table}[!h] 
\caption{The parameters of case studies}
\label{System_para_2}
\centering
\begin{tabular}{clclcl}
\toprule
$m_s$ &$\$$1/(MWh)$^2$ & $\overline {R^{fup}}$& 1.5 MWh & $\Delta t$ & 0.5 hour   \\
$m_g$ &$\$$1/(MWh)$^2$ &$\overline{P^s}$& 3 MW & $N_{OS}$& 150 \\
$m_{gr}$&$\$$15/MWh & $\gamma$ &0.5 & $N_{IT}$& 20  \\
$m_{mgr}$ &$\$$5/MWh & $\epsilon^{jot}$ &0.2 & $N_S$ & 200 \\
\bottomrule
\end{tabular}
\end{table}

The net-load profiles considering solar power generations at load buses are from the UK power networks trial \cite{uk_power_networks_validation_nodate}, and the reserve requirements are from probabilistic forecasts based on these net-load profiles in ref. \cite{ding_additive_2021} with the 90\% prediction intervals. First, we conduct game simulations with different radii of the Wasserstein-metrics ambiguity set $\mathcal P$ and find the optimal radius given the predefined requirement $\epsilon^{jot}$. Then, based on the optimal radius, we simulate games with increasing players ranging from 2 to 14. We assume all the microgrids have the same VoLL of $\$$1/MWh.

%Fig. \ref{net_load_forecast} shows net load forecasts and reserve requirements. Due to solar power and peak loads, reserve requirements are higher from 6 a.m. to 6 p.m. 

%\begin{figure}[!h]
%    \centering
%    \includegraphics[width=3in]{image/network_game.png}
%        \caption{The 30-bus power network game for case studies}
%    \label{IEEE_30}
%\end{figure}

% \begin{figure}[!h]
%    \centering
%    \includegraphics[width=\linewidth]{image/Forecast_plot.png}
%        \caption{(a) An example of probabilistic net load forecast at one bus and (b) total reserve procurement requirements at the reliability levels of 75\% and 90\% }
%    \label{net_load_forecast}
%\end{figure}

%The first case study simulates the game with a fixed number of six microgrids. As shown in Fig. \ref{IEEE_9}, microgrids 1-6 are at buses 21, 23, 24, 26, 29, and 30. Microgrids 1 and 2 have a higher VoLL of $\$$20/MWh, and the rest of the four microgrids have a VoLL of $\$$1/MWh. 

\section{Numerical results} 

The proposed DRJCC game framework is simulated for six hours from 6:00 a.m. to 12:00 a.m.  This end-to-end Bayesian game framework is built in Python using the CVXPY package \cite{diamond_cvxpy_2016} and the scikit-optimize package \cite{noauthor_scikit-optimize_nodate}. The following subsections present the results of case studies.

 \subsection{The optimal radius of Wasserstein-metrics ambiguity sets} 
 
The radius of the Wasserstein-metrics ambiguity determines solution conservativeness, which should be chosen based on grid conditions and system security requirements. We first conduct simulations with an increasing radius ranging from $10^{-5}$ to $10^{-1}$ to find the optimal radius at the given joint violation rate. For comparison, we simulate the benchmark case where there is no system-wise reserve regulation, in other words, relaxing constraints (\ref{DRJCC}) - (\ref{convex_constraint_4}).
  
  \begin{figure}[!h]
    \centering
    \includegraphics[width=\linewidth]{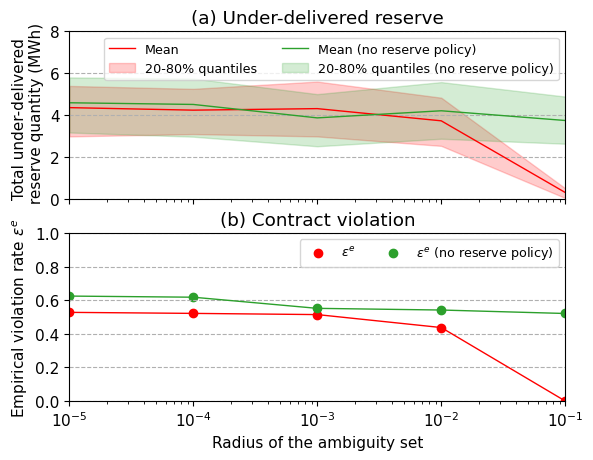}
        \caption{\textcolor{black}{(a) The total amount of under-delivered reserve and (b) the empirical contract violation rates with an increasing radius of the ambiguity set}}
            \label{Reserve_result_0}
\end{figure}

Fig. \ref{Reserve_result_0} shows the under-delivered reserve and joint violation rates from $N_{OS} =150$ out-of-sample tests. The red line is the mean value of the under-delivered reserve from all tests, and the red shade shows the 20-80\% quantile. The proposed game framework is effective when the radius is between $10^{-2}$ and $10^{-1}$, where the game results differ from the benchmark case (i.e., no reserve policy) significantly. We choose an optimal radius of $\nu = 0.035$ for a joint contract violation rate of 0.2.

% \begin{figure}[!h]
%    \centering
%    \includegraphics[width=\linewidth]{image/Profit_results.png}
%        \caption{Average costs of ISO and microgrids with an increasing radius of the ambiguity set}
%    \label{Profit_result_0}
%\end{figure}
%
%Using the mean values in out-of-sample tests, we calculate the average costs of ISO and six microgrids as the objective functions of models M1 and M2, respectively. Fig. \ref{Profit_result_0} summarizes results in five simulations. The costs of microgrids 1 and 2 reduce when the radius increases from $10^{-3}$ to $10^{-1}$, while the costs of the other four microgrids remain constant with the radius increases. This is because microgrids 1 and 2 have the higher VoLL and tend to be more risk-aware (i.e., bid less) when the ambiguity set becomes more conservative (i.e., larger). The average cost of ISO increases rapidly as the ambiguity set becomes larger, as ISO needs more expensive reserve services from thermal generators.

  \subsection{Game convergence with the increasing number of players} 
  
Based on this optimal radius, we simulate games with a varying number of players. For each simulation, we increase two microgrid players (from passive net loads at random buses) until the number of players reaches 14. Since microgrids have the same VoLL, we assume their violation rates are the same (i.e., $\epsilon^{ind}_0$ = $\epsilon^{ind}_1$ =......= $\epsilon^{ind}_n$).

For comparison, we include two scenarios considering the upper and lower bounds of the equation (\ref{range}). The scenario for the lower bound (i.e., $ \boldsymbol \epsilon^{ind} := \frac{\epsilon^{jot}}{N_{MG}} $) denotes C1, and the scenario for the upper bound (i.e., $ \boldsymbol \epsilon^{ind} := \epsilon^{jot}$) denotes C2. The proposed method C3 aims to find the optimal individual violation rate between two values. To mitigate randomness in the outcomes, we choose different initial starting points and repeat simulations five times. 

%We record the minimum value of the evaluation function in each run, as shown in Fig. \ref{time_and_std}. 

\begin{table}[!h]
\caption{The average optimized individual violation allowance over five repeating experiments;}
\label{optimized_rates}
\centering
\begin{tabular}{lccccccc}
\toprule
Players &2& 4& 6 &8& 10 & 12 & 14 \\
\midrule
$ \overline {\boldsymbol \epsilon^{ind}} (\%)$ &20.00 &15.64 &13.76 &11.74& 10.52&  9.74&  9.72\\
\bottomrule
\end{tabular}
\end{table}

%  \begin{figure}[!h]
%    \centering
%    \includegraphics[width=\linewidth]{image/time_and_std.png}
%            \caption{\textcolor{black}{The optimized individual violation rates in simulations with an increasing number of players;  Orange dots show results in each run, and the line shows the average value}}
%    \label{time_and_std}
%\end{figure}

Table \ref{optimized_rates} shows the average optimized individual violation allowance $ \overline {\boldsymbol \epsilon^{ind}}$ for a different number of players over five repeating experiments. The value of $ \boldsymbol \epsilon^{ind}$ decreases with the number of players increasing since the correlation of microgrids' reserve provisions increases. 

% Fig. \ref{time_and_std} shows the distribution of simulation results and average computation time. Two findings can be extracted. First, the value of $ \boldsymbol \epsilon^{ind}$ decreases when the number of players increases since the correlation of microgrids' reserve provisions increases. Second, the computation time increases with the number of players. For 20 iterations, the average computation time is less than 20 minutes. 

  \begin{figure}[!h]
    \centering
    \includegraphics[width=\linewidth]{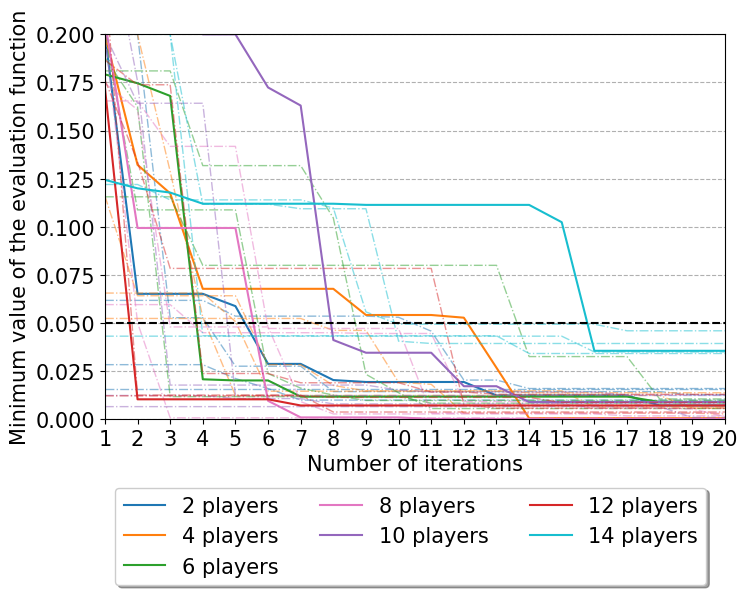}
            \caption{\textcolor{black}{The convergence of Bayesian optimization in games with increasing players. All the games converge to points with a deviation less than 0.05, highlighted in black dashed line; Colored dash lines show the convergence of games in the repeating experiments with different starting points}}
    \label{convergence_plot}
\end{figure}

Fig. \ref{convergence_plot} shows the convergence of Bayesian optimization in different games. The proposed Bayesian optimization allows a maximum of 20 iterations and accepts the game result when the empirical joint violation rate deviates from the predefined requirement by less than 0.05. In all experiments, Bayesian optimization can converge to an optimal point within the maximum acceptable deviation (as shown in the black dashed line in Fig. \ref{convergence_plot}). Given the finite global optimization time, these points are optional to be the global optimum, but they are optimal solutions satisfying both the market clearing time and system reliability requirement. 

  \subsection{Game results from out-of-sample tests}

\begin{figure}[!h]
    \centering
    \includegraphics[width=\linewidth]{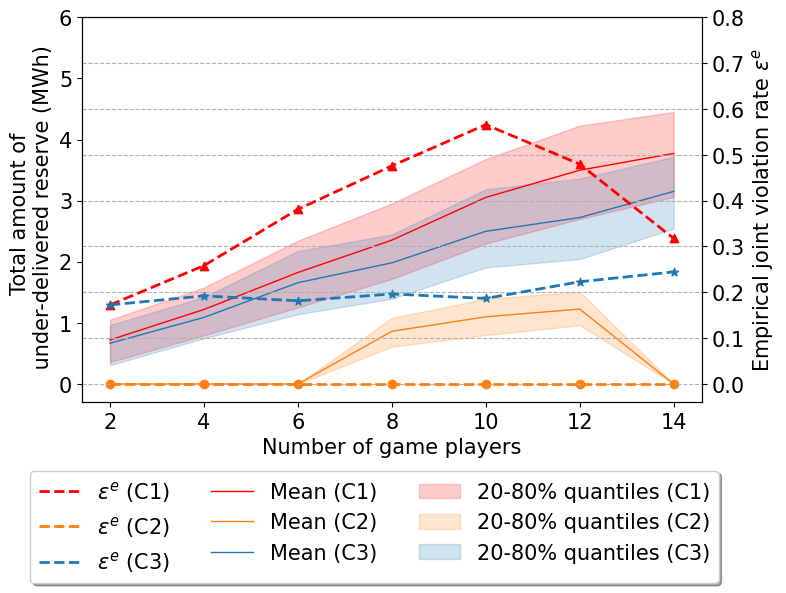}
        \caption{The total amount of under-delivered reserve and empirical joint violation rates in three scenarios}
    \label{Reserve_results_2}
\end{figure}

Fig. \ref{Reserve_results_2} shows how the under-delivered reserve and joint violation rate change in three scenarios. Only the proposed method can regulate the under-delivered reserve and secure the joint violation rate of all microgrids in a game. Using the Bayesian optimization C3 (blue), the total under-delivered reserve increases steadily, and the empirical joint violation rate is regulated at $0.2 \pm 0.05$.

In scenario C1, the total amount of under-delivered reserve increases rapidly, and the joint violation rate exceeds 0.2 even for the 4-player game and quickly increases to 0.6 in the 10-player game. This means that although C1 using single CCs could increase microgrids' profits and lower ISO costs, the result is unreliable for grid operations. In scenario C2, on the contrary, the joint violation rate stays close to zero. Neglecting the correlation leads to an overly conservative reserve regulation and almost no reserve bids in scenario C3. 

  \begin{figure}[!h]
    \centering
    \includegraphics[width=\linewidth]{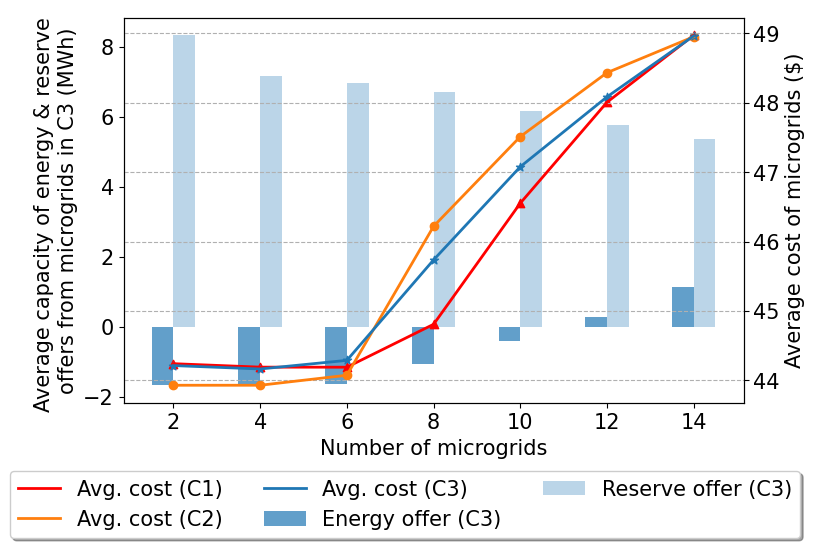}
        \caption{Average costs of microgrids in three scenarios and their average capacity of energy and reserve offers in C3 }
    \label{players_results}
\end{figure}

Fig. \ref{players_results} summarizes the average cost of microgrids and the average capacity of their energy and reserve offers in scenario C3. The positive value represents that microgrids provide services to the grid and vice versa. There is a trade-off between energy and reserve services, and the reserve capacity of each microgrid decreases as the number of microgrids increases due to market competition.

  \begin{figure}[!h]
    \centering
    \includegraphics[width=\linewidth]{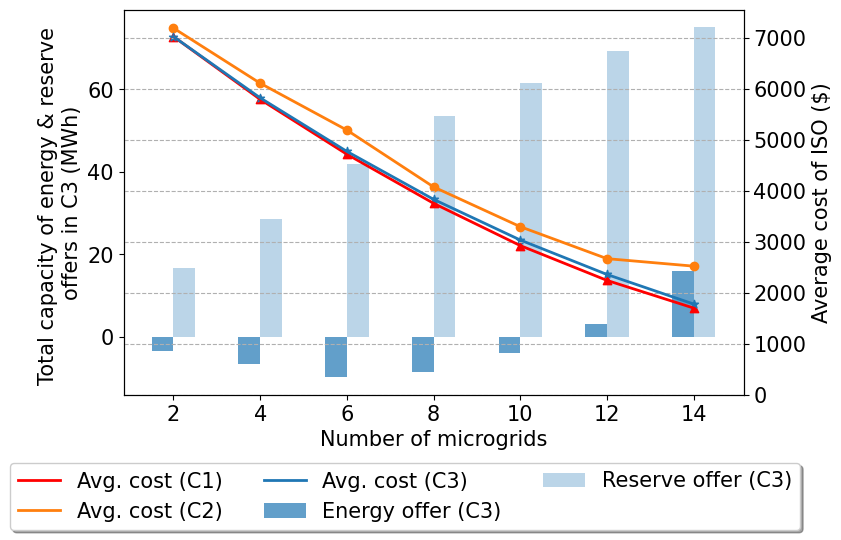}
        \caption{Average costs of ISO in three scenarios and the total capacity of energy and reserve offers in C3}
    \label{iso_results}
\end{figure}

The game starts with two microgrid players. At first, the average cost of microgrids in scenarios C1 and C3 is more than in C2. As in scenarios C1 and C3, microgrids prefer to invest in the reserve market for higher profits, even if they import power from the grid (as shown in light and dark blue bars). They gain a high load-shedding penalty in return. As the players increase, microgrids reduce their reserve offers due to market competition. The average cost of microgrids in scenarios C1 and C3 becomes lower than in C2. If the number increases to 10 or more, the reserve market becomes saturated, and microgrids need to sell energy for more profits.
 
Fig. \ref{iso_results} summarizes the cost of ISO and the total energy and reserve bidding capacity from microgrids in scenario C3. The total reserve capacity of microgrids increases, replacing the expensive reserve services from thermal generators (as presented in the light blue bars). Therefore, the cost of ISO directly decreases with the number of microgrid players increasing in all three scenarios. 

%Three scenarios ranked by the ISO's cost are C2, C3, and C1, which reflect the conservativeness of system reserve regulations using three different approaches.

  \subsection{Comparision of Bayesian optimization and grid search}

% \begin{table}[!h]
%\caption{The optimized individual violation rate by grid search}
%\label{optimized_rates_gs}
%\centering
%\begin{tabular}{lccccccc}
%\toprule
%Players &2& 4& 6 &8& 10 & 12 & 14 \\
%\midrule
%$ \boldsymbol \epsilon^{ind'} (\%)$ &20.00 &15.20 &14.53 &11.28& 11.00&  10.82&  9.79\\
%\bottomrule
%\end{tabular}
%\end{table}
%  

We compare Bayesian optimization with the naive grid search method. Grid search partitions the exploration space into a grid and exhaustively searches every grid point. The explore space (\ref{range}) is partitioned into $N_{sch}=20$ discrete search points with an equal distance $d= \frac{(N_{MG}-1)\epsilon^{jot}}{N_{MG}(N_{sch}-1)}$. We then solve model M3 at each search point $ \epsilon_m \in \{\frac{\epsilon^{jot}}{N_{MG}}, \frac{\epsilon^{jot}}{N_{MG}} + d ,\cdots, \epsilon^{jot}\}$ and identify the optimal violation rate.

  \begin{figure}[!h]
    \centering
    \includegraphics[width=\linewidth]{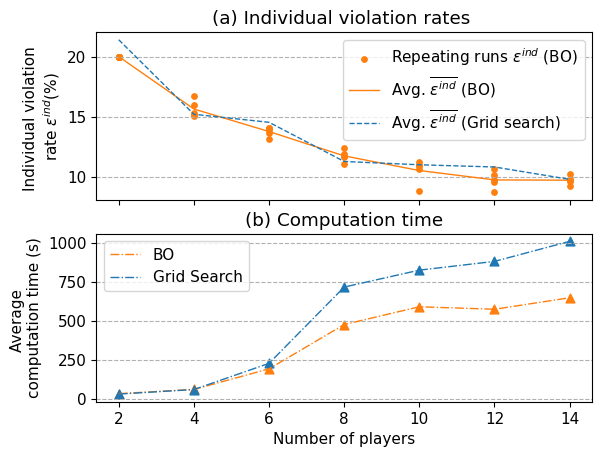}
        \caption{(a) The optimal individual violation rates and (b) the computation time of the game simulations using Bayesian optimization (blue) and grid search (orange)}
    \label{comparision_time}
\end{figure}

 The grid search obtains similar game results to Bayesian optimization in all simulations, as shown in Fig. \ref{comparision_time} (a). However, grid search takes longer in solution convergence compared with Bayesian optimization, and this time difference increases with the number of players (Fig. \ref{comparision_time} (b)). For the 14-player game, Bayesian optimization is faster than the grid search by 54\% (i.e., 361 seconds). This is because the grid search exhaustively searches the possibility and is less efficient than Bayesian optimization.

  \section{Conclusions} 
  
We propose a DRJCC Stackelberg game to coordinate renewable microgrids in energy and reserve markets, considering service delivery risks. In light of reliable reserve provision considering solar power uncertainty, a two-stage DRJCC based on the Wasserstein-metric ambiguity set models the worst-case reserve performance and the expected penalty of each microgrid due to undelivered services. 

This DRJCC game problem is intractable when the individual contract violation rate is a design variable. Considering the unspecified and changing correlation in games, these optimal individual violation rates should precisely regulate the under-delivered reserve and ensure the profits of microgrids. Instead of enforcing conservative assumptions on the DRJCC problem in the literature, we design a novel Bayesian optimization method to approach the optimal solutions based on the existence and uniqueness of game equilibrium. We include two extreme scenarios with a prior value of individual violation rates to compare with the proposed Bayesian optimization (C3). These two scenarios assume players have the same performances (C1) or no correlation in reserve provision performances (C2). 

Games are simulated with different numbers of players in the repeating experiments, which range from 2 to 14. Only the Bayesian optimization approach (C3) can converge to the optimal individual violation rates so that the empirical joint violation rates only deviate from the predefined requirement by $0.05$ in 20 iterations. In contrast, the other two approaches (C1 and C2) are either highly unreliable for the grid or too conservative to generate reserve bids. 

Bayesian optimization is then compared to the grid search method. It can find global optimality faster than the grid search by leveraging prior knowledge during exploration, especially when the problem size increases. Nevertheless, due to the rapidly increasing computation costs for high-dimension data, it is still worth considering dimension reduction to reduce the consequential computation costs when using Bayesian optimization. 
 
 \section{Appendix A: Proofs for the existence and uniqueness of game solution} 

We use the variational inequality framework \cite{gabriel_equilibria_2013} $VI(X, F)$ to prove the existence and uniqueness of the game solution. This is given by,

    \begin{equation} 
F( \boldsymbol x^{*})^T ( \boldsymbol x  -  \boldsymbol x^{*}) \geq 0 \quad \text{for all} \quad  \boldsymbol x \in \mathcal{K}
\label{upper}
 \end{equation}

where the strategy set $\mathcal{K}$ and the game map $F( \boldsymbol x)$ when the proposed game has $n$ microgrid players are given by,

 \begin{equation} 
\begin{aligned}
\mathcal K & =  [\boldsymbol x^{H}_1, \boldsymbol x^{L}_1 ,.....,  \boldsymbol x^{L}_n ] \\
& =  [P^g, R^g, \pi^r, R^{up},  \pi^e,  P^{ex}, P^s_1, \pi^r_1, R^{up}_1,  \cdots,  \textbf W^l_n ]
\end{aligned}
\label{discrete_points}
 \end{equation}

\begin{equation} 
F( \boldsymbol{x}) = \begin{pmatrix}
\nabla \mathcal{J}^{H}_1 (\boldsymbol x^{H}_1, \boldsymbol x^{H}_{-1}) \\
\nabla \mathcal{J}^{L}_1 (\boldsymbol x^{L}_1, \boldsymbol x^{L}_{-1}) \\
... ...\\
\nabla \mathcal{J}^{L}_n (\boldsymbol x^{L}_n, \boldsymbol x^{L}_{-n}) \end{pmatrix} = \begin{pmatrix}
\sum 2m_k^g P^g t \\
\sum  m_k^{gr} t \\
... ...\\
% \nu\\
\frac{1}{N_S} \end{pmatrix}
\label{game_map}
 \end{equation}
 
The game map $F( \boldsymbol{x})$ is continuously differentiable since all players' cost functions $\mathcal{J}$ are quadratic or linear functions. In addition, the strategy set $\mathcal K$ is compact, convex, and non-empty. Therefore, a solution set $\mathcal{S} (\mathcal K, F)$ exists for the game. We then derive the game Jacobian matrix of $\nabla F( \boldsymbol x)$ to observe the properties of Nash equilibrium.

\begin{equation}
\footnotesize
\begin{aligned}
&\nabla F( \boldsymbol x) = \\
&\begin{bmatrix}
\textcolor{blue}{\sum 2m_k^g  t} &0 & 0 & 0 & 0 & 0 & 0 & 0 &\cdots & 0\\
0 & \textcolor{blue}{0} & 0 & 0 & 0& 0 & 0 & 0 & \cdots & 0\\
0 & 0 & \textcolor{blue}{0} & \sum t  & 0 & 0&  0 & 0 & \cdots & 0\\
0 & 0 & \sum t  & \textcolor{blue}{0}  &  0   & 0 & 0 & 0 & \cdots & 0\\
0 & 0 & 0 & 0 & \textcolor{blue}{0} & -\sum t & 0 & 0 & \cdots & 0\\
0 & 0 & 0 & 0  & -\sum t  & \textcolor{blue}{0} &  0 & 0 & \cdots & 0\\
0 & 0 & 0 & 0  & 0 & 0 &  \textcolor{blue}{0}  & 0 & \cdots & 0\\
0 & 0 & 0 & 0 & 0 & 0 &  0 &  \textcolor{blue}{\sum 2m_i^s  t} & \cdots & 0\\
\vdots & \vdots & \vdots & \vdots &  \vdots &  \vdots &  \vdots &  \vdots & \ddots & \vdots\\
0 & 0 & 0 & 0 & 0 & 0 &  0 &  0 & \cdots &  \textcolor{blue}{0} \\
\end{bmatrix}
\end{aligned}
\end{equation}

The Jacobian matrix above is symmetric, indicated by the blue diagonal entries. This means that the corresponding game is \textit{integrable} \cite{gabriel_equilibria_2013}. Based on the Principle of Integrability, setting $x^0 = 0$, the equivalent optimization problem has an objective function given by,

\begin{equation} 
\begin{aligned}
\theta (\boldsymbol x) = &\int_{1}^{0}  F( x^0 + \gamma (x-x^0))^{\top} (x-x^0) d\gamma\\
= & \sum_{t\in \mathcal T}\sum_{k\in \mathcal K}  m_{k}^g ( \Delta t P^g_{t,k})^2 +   \frac{1}{2} \Delta t m_{k}^{gr}R^g_{t, k} \\
&+ \sum_{t\in \mathcal T} \sum_{i\in \mathcal I} m_i^s (\Delta t P_{t, i}^s) ^2 +   \frac{1}{2} \Delta t  m_i^{mgr} R_{t, i}^{up} \\
& + \sum_{t\in \mathcal T} \sum_{i\in \mathcal I} (\nu l^l_i + \frac{\boldsymbol 1^{\top}\textbf W^l_i}{N_S})
 \end{aligned}
  \label{objective3}
 \end{equation}
 
  \qquad \textrm{s.t. \quad equations} \quad (\ref{power_balance}) - (\ref{convex_constraint_4}) \quad \textrm{and} \quad (\ref{lower_0}) - (\ref{solar_reserve_bid})
 \\

The objective function is quadratic. Given a fixed value of $\epsilon^{ind}_i$ in each iteration of Bayesian optimization, the equivalent optimization problem is convex. Therefore, the $VI(X, F)$ admits a unique solution.

 \bibliographystyle{IEEEtran}
\bibliography{bibtex/IEEEconference}

% Generated by IEEEtran.bst, version: 1.14 (2015/08/26)
\begin{thebibliography}{10}
\providecommand{\url}[1]{#1}
\csname url@samestyle\endcsname
\providecommand{\newblock}{\relax}
\providecommand{\bibinfo}[2]{#2}
\providecommand{\BIBentrySTDinterwordspacing}{\spaceskip=0pt\relax}
\providecommand{\BIBentryALTinterwordstretchfactor}{4}
\providecommand{\BIBentryALTinterwordspacing}{\spaceskip=\fontdimen2\font plus
\BIBentryALTinterwordstretchfactor\fontdimen3\font minus
  \fontdimen4\font\relax}
\providecommand{\BIBforeignlanguage}[2]{{%
\expandafter\ifx\csname l@#1\endcsname\relax
\typeout{** WARNING: IEEEtran.bst: No hyphenation pattern has been}%
\typeout{** loaded for the language `#1'. Using the pattern for}%
\typeout{** the default language instead.}%
\else
\language=\csname l@#1\endcsname
\fi
#2}}
\providecommand{\BIBdecl}{\relax}
\BIBdecl

\bibitem{department_for_business_energy__industrial_strategy_powering_nodate}
\BIBentryALTinterwordspacing
{Department for Business, Energy \& Industrial Strategy}, ``Powering our {Net}
  {Zero} {Future},'' Tech. Rep. [Online]. Available:
  \url{https://www.gov.uk/government/publications/energy-white-paper-powering-our-net-zero-future}
\BIBentrySTDinterwordspacing

\bibitem{keay_electricity_2013}
\BIBentryALTinterwordspacing
M.~Keay, J.~Rhys, and D.~Robinson, ``\BIBforeignlanguage{en}{Electricity
  {Market} {Reform} in {Britain}},'' in \emph{\BIBforeignlanguage{en}{Evolution
  of {Global} {Electricity} {Markets}}}.\hskip 1em plus 0.5em minus 0.4em\relax
  Elsevier, 2013, pp. 31--57. [Online]. Available:
  \url{https://linkinghub.elsevier.com/retrieve/pii/B978012397891200002X}
\BIBentrySTDinterwordspacing

\bibitem{gabriel_equilibria_2013}
\BIBentryALTinterwordspacing
S.~A. Gabriel, A.~J. Conejo, J.~D. Fuller, B.~F. Hobbs, and C.~Ruiz,
  ``Equilibria and {Complementarity} {Problems},'' in \emph{Complementarity
  {Modeling} in {Energy} {Markets}}.\hskip 1em plus 0.5em minus 0.4em\relax New
  York, NY: Springer New York, 2013, vol. 180, pp. 127--179, series Title:
  International Series in Operations Research \& Management Science. [Online].
  Available: \url{http://link.springer.com/10.1007/978-1-4419-6123-5_4}
\BIBentrySTDinterwordspacing

\bibitem{delatorre_finding_2004}
\BIBentryALTinterwordspacing
S.~delaTorre, J.~Contreras, and A.~Conejo, ``\BIBforeignlanguage{en}{Finding
  {Multiperiod} {Nash} {Equilibria} in {Pool}-{Based} {Electricity}
  {Markets}},'' \emph{\BIBforeignlanguage{en}{IEEE Transactions on Power
  Systems}}, vol.~19, no.~1, pp. 643--651, Feb. 2004. [Online]. Available:
  \url{http://ieeexplore.ieee.org/document/1266624/}
\BIBentrySTDinterwordspacing

\bibitem{li_data-driven_2022}
\BIBentryALTinterwordspacing
X.~Li, C.~Li, G.~Chen, and Z.~Y. Dong, ``A {Data}-driven {Joint}
  {Chance}-constrained {Game} for {Renewable} {Energy} {Aggregators} in the
  {Local} {Market},'' \emph{IEEE Transactions on Smart Grid}, pp. 1--1, 2022.
  [Online]. Available: \url{https://ieeexplore.ieee.org/document/9744121/}
\BIBentrySTDinterwordspacing

\bibitem{saad_game-theoretic_2012}
\BIBentryALTinterwordspacing
W.~Saad, Z.~Han, H.~Poor, and T.~Basar, ``Game-{Theoretic} {Methods} for the
  {Smart} {Grid}: {An} {Overview} of {Microgrid} {Systems}, {Demand}-{Side}
  {Management}, and {Smart} {Grid} {Communications},'' \emph{IEEE Signal
  Processing Magazine}, vol.~29, no.~5, pp. 86--105, Sep. 2012. [Online].
  Available: \url{http://ieeexplore.ieee.org/document/6279592/}
\BIBentrySTDinterwordspacing

\bibitem{saguan_market_2006}
\BIBentryALTinterwordspacing
M.~Saguan, N.~Keseric, P.~Dessante, and J.-m. Glachant, ``Market {Power} in
  {Power} {Markets}: {Game} {Theory} vs. {Agent}-{Based} {Approach},'' in
  \emph{2006 {IEEE}/{PES} {Transmission} \& {Distribution} {Conference} and
  {Exposition}: {Latin} {America}}.\hskip 1em plus 0.5em minus 0.4em\relax
  Caracas, Venezuela: IEEE, 2006, pp. 1--6. [Online]. Available:
  \url{http://ieeexplore.ieee.org/document/4104670/}
\BIBentrySTDinterwordspacing

\bibitem{jiang_flexibility_2022}
\BIBentryALTinterwordspacing
T.~Jiang, C.~Wu, R.~Zhang, X.~Li, H.~Chen, and G.~Li, ``Flexibility {Clearing}
  in {Joint} {Energy} and {Flexibility} {Markets} {Considering} {TSO}-{DSO}
  {Coordination},'' \emph{IEEE Transactions on Smart Grid}, pp. 1--1, 2022.
  [Online]. Available: \url{https://ieeexplore.ieee.org/document/9722366/}
\BIBentrySTDinterwordspacing

\bibitem{hobbs_strategic_2000-1}
\BIBentryALTinterwordspacing
B.~Hobbs, C.~Metzler, and J.-S. Pang, ``Strategic gaming analysis for electric
  power systems: an {MPEC} approach,'' \emph{IEEE Transactions on Power
  Systems}, vol.~15, no.~2, pp. 638--645, May 2000. [Online]. Available:
  \url{http://ieeexplore.ieee.org/document/867153/}
\BIBentrySTDinterwordspacing

\bibitem{li_distributed_2023}
\BIBentryALTinterwordspacing
Z.~Li, L.~Wu, Y.~Xu, L.~Wang, and N.~Yang,
  ``\BIBforeignlanguage{en}{Distributed tri-layer risk-averse stochastic game
  approach for energy trading among multi-energy microgrids},''
  \emph{\BIBforeignlanguage{en}{Applied Energy}}, vol. 331, p. 120282, Feb.
  2023. [Online]. Available:
  \url{https://linkinghub.elsevier.com/retrieve/pii/S0306261922015392}
\BIBentrySTDinterwordspacing

\bibitem{bruninx_interaction_2020}
\BIBentryALTinterwordspacing
K.~Bruninx, H.~Pandzic, H.~Le~Cadre, and E.~Delarue, ``On the {Interaction}
  {Between} {Aggregators}, {Electricity} {Markets} and {Residential} {Demand}
  {Response} {Providers},'' \emph{IEEE Transactions on Power Systems}, vol.~35,
  no.~2, pp. 840--853, Mar. 2020. [Online]. Available:
  \url{https://ieeexplore.ieee.org/document/8848616/}
\BIBentrySTDinterwordspacing

\bibitem{li_risk-averse_2017}
\BIBentryALTinterwordspacing
C.~Li, Y.~Xu, X.~Yu, C.~Ryan, and T.~Huang, ``Risk-{Averse} {Energy} {Trading}
  in {Multienergy} {Microgrids}: {A} {Two}-{Stage} {Stochastic} {Game}
  {Approach},'' \emph{IEEE Transactions on Industrial Informatics}, vol.~13,
  no.~5, pp. 2620--2630, Oct. 2017. [Online]. Available:
  \url{http://ieeexplore.ieee.org/document/8010361/}
\BIBentrySTDinterwordspacing

\bibitem{vespermann_risk_2021}
\BIBentryALTinterwordspacing
N.~Vespermann, T.~Hamacher, and J.~Kazempour, ``Risk {Trading} in {Energy}
  {Communities},'' \emph{IEEE Transactions on Smart Grid}, vol.~12, no.~2, pp.
  1249--1263, Mar. 2021. [Online]. Available:
  \url{https://ieeexplore.ieee.org/document/9222043/}
\BIBentrySTDinterwordspacing

\bibitem{ordoudis_energy_2021}
\BIBentryALTinterwordspacing
C.~Ordoudis, V.~A. Nguyen, D.~Kuhn, and P.~Pinson,
  ``\BIBforeignlanguage{en}{Energy and reserve dispatch with distributionally
  robust joint chance constraints},'' \emph{\BIBforeignlanguage{en}{Operations
  Research Letters}}, vol.~49, no.~3, pp. 291--299, May 2021. [Online].
  Available:
  \url{https://linkinghub.elsevier.com/retrieve/pii/S0167637721000213}
\BIBentrySTDinterwordspacing

\bibitem{ding_distributionally_2022}
\BIBentryALTinterwordspacing
Y.~Ding, T.~Morstyn, and M.~D. McCulloch, ``Distributionally {Robust} {Joint}
  {Chance}-{Constrained} {Optimization} for {Networked} {Microgrids}
  {Considering} {Contingencies} and {Renewable} {Uncertainty},'' \emph{IEEE
  Transactions on Smart Grid}, vol.~13, no.~3, pp. 2467--2478, May 2022.
  [Online]. Available: \url{https://ieeexplore.ieee.org/document/9709590/}
\BIBentrySTDinterwordspacing

\bibitem{mazadi_impact_2013}
\BIBentryALTinterwordspacing
M.~Mazadi, W.~D. Rosehart, H.~Zareipour, O.~P. Malik, and M.~Oloomi,
  ``\BIBforeignlanguage{en}{Impact of wind integration on electricity markets:
  a chance-constrained {Nash} {Cournot} model: {CHANCE}-{CONSTRAINED} {NASH}
  {COURNOT} {MODEL}},'' \emph{\BIBforeignlanguage{en}{International
  Transactions on Electrical Energy Systems}}, vol.~23, no.~1, pp. 83--96, Jan.
  2013. [Online]. Available:
  \url{https://onlinelibrary.wiley.com/doi/10.1002/etep.650}
\BIBentrySTDinterwordspacing

\bibitem{xie_optimized_2022}
\BIBentryALTinterwordspacing
W.~Xie, S.~Ahmed, and R.~Jiang, ``\BIBforeignlanguage{en}{Optimized
  {Bonferroni} approximations of distributionally robust joint chance
  constraints},'' \emph{\BIBforeignlanguage{en}{Mathematical Programming}},
  vol. 191, no.~1, pp. 79--112, Jan. 2022. [Online]. Available:
  \url{https://link.springer.com/10.1007/s10107-019-01442-8}
\BIBentrySTDinterwordspacing

\bibitem{chen_cvar_2010}
\BIBentryALTinterwordspacing
W.~Chen, M.~Sim, J.~Sun, and C.-P. Teo, ``\BIBforeignlanguage{en}{From {CVaR}
  to {Uncertainty} {Set}: {Implications} in {Joint} {Chance}-{Constrained}
  {Optimization}},'' \emph{\BIBforeignlanguage{en}{Operations Research}},
  vol.~58, no.~2, pp. 470--485, Apr. 2010. [Online]. Available:
  \url{http://pubsonline.informs.org/doi/10.1287/opre.1090.0712}
\BIBentrySTDinterwordspacing

\bibitem{baker_joint_2019}
\BIBentryALTinterwordspacing
K.~Baker and A.~Bernstein, ``Joint {Chance} {Constraints} in {AC} {Optimal}
  {Power} {Flow}: {Improving} {Bounds} {Through} {Learning},'' \emph{IEEE
  Transactions on Smart Grid}, vol.~10, no.~6, pp. 6376--6385, Nov. 2019.
  [Online]. Available: \url{https://ieeexplore.ieee.org/document/8662704/}
\BIBentrySTDinterwordspacing

\bibitem{jia_iterative_2021}
M.~Jia, G.~Hug, and C.~Shen, ``Iterative {Decomposition} of {Joint} {Chance}
  {Constraints} in {OPF},'' \emph{IEEE Transactions on Power Systems}, vol.~36,
  no.~5, pp. 4836--4839, Sep. 2021.

\bibitem{arrigo_embedding_2022}
\BIBentryALTinterwordspacing
A.~Arrigo, J.~Kazempour, Z.~De~Greve, J.-F. Toubeau, and F.~Vallee, ``Embedding
  {Dependencies} {Between} {Wind} {Farms} in {Distributionally} {Robust}
  {Optimal} {Power} {Flow},'' \emph{IEEE Transactions on Power Systems}, pp.
  1--14, 2022. [Online]. Available:
  \url{https://ieeexplore.ieee.org/document/9944957/}
\BIBentrySTDinterwordspacing

\bibitem{inatsu2022bayesian}
Y.~Inatsu, S.~Takeno, M.~Karasuyama, and I.~Takeuchi, ``Bayesian optimization
  for distributionally robust chance-constrained problem,'' in
  \emph{International Conference on Machine Learning}.\hskip 1em plus 0.5em
  minus 0.4em\relax PMLR, 2022, pp. 9602--9621.

\bibitem{shapiro_bayesian_2023}
\BIBentryALTinterwordspacing
A.~Shapiro, E.~Zhou, and Y.~Lin, ``\BIBforeignlanguage{en}{Bayesian
  {Distributionally} {Robust} {Optimization}},''
  \emph{\BIBforeignlanguage{en}{SIAM Journal on Optimization}}, vol.~33, no.~2,
  pp. 1279--1304, Jun. 2023. [Online]. Available:
  \url{https://epubs.siam.org/doi/10.1137/21M1465548}
\BIBentrySTDinterwordspacing

\bibitem{ding_coordinating_2023}
\BIBentryALTinterwordspacing
Y.~Ding, S.~Wang, and B.~Hobbs, ``\BIBforeignlanguage{en}{Coordinating
  renewable microgrids for reliable reserve services: a distributionally robust
  chance-constrained game},'' in \emph{\BIBforeignlanguage{en}{Proceedings of
  the 14th {ACM} {International} {Conference} on {Future} {Energy}
  {Systems}}}.\hskip 1em plus 0.5em minus 0.4em\relax Orlando FL USA: ACM, Jun.
  2023, pp. 324--332. [Online]. Available:
  \url{https://dl.acm.org/doi/10.1145/3575813.3597342}
\BIBentrySTDinterwordspacing

\bibitem{nemirovski_convex_2007}
\BIBentryALTinterwordspacing
A.~Nemirovski and A.~Shapiro, ``\BIBforeignlanguage{en}{Convex {Approximations}
  of {Chance} {Constrained} {Programs}},'' \emph{\BIBforeignlanguage{en}{SIAM
  Journal on Optimization}}, vol.~17, no.~4, pp. 969--996, Jan. 2007. [Online].
  Available: \url{http://epubs.siam.org/doi/10.1137/050622328}
\BIBentrySTDinterwordspacing

\bibitem{ruiz_tutorial_2014}
\BIBentryALTinterwordspacing
C.~Ruiz, A.~J. Conejo, J.~D. Fuller, S.~A. Gabriel, and B.~F. Hobbs,
  ``\BIBforeignlanguage{en}{A tutorial review of complementarity models for
  decision-making in energy markets},'' \emph{\BIBforeignlanguage{en}{EURO
  Journal on Decision Processes}}, vol.~2, no. 1-2, pp. 91--120, Jun. 2014.
  [Online]. Available:
  \url{https://linkinghub.elsevier.com/retrieve/pii/S2193943821000285}
\BIBentrySTDinterwordspacing

\bibitem{gao_distributionally_2023}
\BIBentryALTinterwordspacing
R.~Gao and A.~Kleywegt, ``\BIBforeignlanguage{en}{Distributionally {Robust}
  {Stochastic} {Optimization} with {Wasserstein} {Distance}},''
  \emph{\BIBforeignlanguage{en}{Mathematics of Operations Research}}, vol.~48,
  no.~2, pp. 603--655, May 2023. [Online]. Available:
  \url{https://pubsonline.informs.org/doi/10.1287/moor.2022.1275}
\BIBentrySTDinterwordspacing

\bibitem{lovric_kullback-leibler_2011}
\BIBentryALTinterwordspacing
J.~M. Joyce, ``\BIBforeignlanguage{en}{Kullback-{Leibler} {Divergence}},'' in
  \emph{\BIBforeignlanguage{en}{International {Encyclopedia} of {Statistical}
  {Science}}}, M.~Lovric, Ed.\hskip 1em plus 0.5em minus 0.4em\relax Berlin,
  Heidelberg: Springer Berlin Heidelberg, 2011, pp. 720--722. [Online].
  Available: \url{http://link.springer.com/10.1007/978-3-642-04898-2_327}
\BIBentrySTDinterwordspacing

\bibitem{netessine_wasserstein_2019}
\BIBentryALTinterwordspacing
D.~Kuhn, P.~M. Esfahani, V.~A. Nguyen, and S.~Shafieezadeh-Abadeh,
  ``\BIBforeignlanguage{en}{Wasserstein {Distributionally} {Robust}
  {Optimization}: {Theory} and {Applications} in {Machine} {Learning}},'' in
  \emph{\BIBforeignlanguage{en}{Operations {Research} \& {Management} {Science}
  in the {Age} of {Analytics}}}, S.~Netessine, D.~Shier, and H.~J. Greenberg,
  Eds.\hskip 1em plus 0.5em minus 0.4em\relax INFORMS, Oct. 2019, pp. 130--166.
  [Online]. Available:
  \url{http://pubsonline.informs.org/doi/10.1287/educ.2019.0198}
\BIBentrySTDinterwordspacing

\bibitem{california_iso_oasis_nodate}
\BIBentryALTinterwordspacing
{California ISO}, ``{OASIS} - {OASIS} {Prod} - {PUBLIC} - 0.'' [Online].
  Available: \url{http://oasis.caiso.com/mrioasis/logon.do}
\BIBentrySTDinterwordspacing

\bibitem{zymler_distributionally_2013}
\BIBentryALTinterwordspacing
S.~Zymler, D.~Kuhn, and B.~Rustem, ``\BIBforeignlanguage{en}{Distributionally
  robust joint chance constraints with second-order moment information},''
  \emph{\BIBforeignlanguage{en}{Mathematical Programming}}, vol. 137, no. 1-2,
  pp. 167--198, Feb. 2013. [Online]. Available:
  \url{http://link.springer.com/10.1007/s10107-011-0494-7}
\BIBentrySTDinterwordspacing

\bibitem{mohajerin_esfahani_data-driven_2018}
\BIBentryALTinterwordspacing
P.~Mohajerin~Esfahani and D.~Kuhn, ``\BIBforeignlanguage{en}{Data-driven
  distributionally robust optimization using the {Wasserstein} metric:
  performance guarantees and tractable reformulations},''
  \emph{\BIBforeignlanguage{en}{Mathematical Programming}}, vol. 171, no. 1-2,
  pp. 115--166, Sep. 2018. [Online]. Available:
  \url{http://link.springer.com/10.1007/s10107-017-1172-1}
\BIBentrySTDinterwordspacing

\bibitem{arrigo_wasserstein_2022}
\BIBentryALTinterwordspacing
A.~Arrigo, C.~Ordoudis, J.~Kazempour, Z.~De~Grève, J.-F. Toubeau, and
  F.~Vallée, ``\BIBforeignlanguage{en}{Wasserstein distributionally robust
  chance-constrained optimization for energy and reserve dispatch: {An} exact
  and physically-bounded formulation},'' \emph{\BIBforeignlanguage{en}{European
  Journal of Operational Research}}, vol. 296, no.~1, pp. 304--322, Jan. 2022.
  [Online]. Available:
  \url{https://linkinghub.elsevier.com/retrieve/pii/S0377221721003271}
\BIBentrySTDinterwordspacing

\bibitem{masaki_adachi_sober_nodate}
\BIBentryALTinterwordspacing
{Masaki Adachi}, {Satoshi Hayakawa}, {Saad Hamid}, {Martin Jørgensen}, {Harald
  Oberhauser}, and {Micheal A. Osborne}, ``{SOBER}: {Scalable} {Batch}
  {Bayesian} {Optimization} and {Quadrature} using {Recombination}
  {Constraints}.'' [Online]. Available: \url{https://arxiv.org/abs/2301.11832}
\BIBentrySTDinterwordspacing

\bibitem{shahriari_taking_2016}
\BIBentryALTinterwordspacing
B.~Shahriari, K.~Swersky, Z.~Wang, R.~P. Adams, and N.~de~Freitas, ``Taking the
  {Human} {Out} of the {Loop}: {A} {Review} of {Bayesian} {Optimization},''
  \emph{Proceedings of the IEEE}, vol. 104, no.~1, pp. 148--175, Jan. 2016.
  [Online]. Available: \url{https://ieeexplore.ieee.org/document/7352306/}
\BIBentrySTDinterwordspacing

\bibitem{wilson_maximizing_2018}
\BIBentryALTinterwordspacing
J.~Wilson, F.~Hutter, and M.~Deisenroth, ``Maximizing acquisition functions for
  {Bayesian} optimization,'' in \emph{Advances in {Neural} {Information}
  {Processing} {Systems}}, S.~Bengio, H.~Wallach, H.~Larochelle, K.~Grauman,
  N.~Cesa-Bianchi, and R.~Garnett, Eds., vol.~31.\hskip 1em plus 0.5em minus
  0.4em\relax Curran Associates, Inc., 2018. [Online]. Available:
  \url{https://proceedings.neurips.cc/paper/2018/file/498f2c21688f6451d9f5fd09d53edda7-Paper.pdf}
\BIBentrySTDinterwordspacing

\bibitem{power_systems_test_case_archive_30_nodate}
\BIBentryALTinterwordspacing
{Power Systems Test Case Archive}, ``30 {Bus} {Power} {Flow} {Test} {Case}.''
  [Online]. Available: \url{http://labs.ece.uw.edu/pstca/pf30/pg_tca30bus.htm}
\BIBentrySTDinterwordspacing

\bibitem{uk_power_networks_validation_nodate}
\BIBentryALTinterwordspacing
{UK Power Networks}, ``Validation of {Photovoltaic} ({PV}) {Connection}
  {Assessment} {Tool},'' Tech. Rep. [Online]. Available:
  \url{https://www.ofgem.gov.uk/ofgem-publications/93938/ pvtoolcdrfinal-pdf}
\BIBentrySTDinterwordspacing

\bibitem{ding_additive_2021}
\BIBentryALTinterwordspacing
Y.~Ding and M.~McCulloch, ``\BIBforeignlanguage{en}{Additive {Gaussian} process
  prediction for electrical loads compared with deep learning models},'' in
  \emph{\BIBforeignlanguage{en}{Proceedings of the {Twelfth} {ACM}
  {International} {Conference} on {Future} {Energy} {Systems}}}.\hskip 1em plus
  0.5em minus 0.4em\relax Virtual Event Italy: ACM, Jun. 2021, pp. 499--506.
  [Online]. Available: \url{https://dl.acm.org/doi/10.1145/3447555.3466592}
\BIBentrySTDinterwordspacing

\bibitem{diamond_cvxpy_2016}
S.~Diamond and S.~Boyd, ``{CVXPY}: a python-embedded modeling language for
  convex optimization,'' \emph{The Journal of Machine Learning Research},
  vol.~17, no.~1, pp. 2909--2913, Jan. 2016.

\bibitem{noauthor_scikit-optimize_nodate}
\BIBentryALTinterwordspacing
``scikit-optimize.'' [Online]. Available:
  \url{https://scikit-optimize.github.io/stable/}
\BIBentrySTDinterwordspacing

\end{thebibliography}

% that's all, folks
\end{document}